\documentclass{article}

\usepackage{graphicx}
\usepackage{amsmath, amsthm}
\usepackage{amsfonts}
\usepackage{amssymb}
\usepackage{url}
\usepackage{hyperref}  
\usepackage{verbatim}
\usepackage{color}
\definecolor{red}{rgb}{1,0,0}

\definecolor{blu}{rgb}{0,0,1}

\usepackage[mathlines]{lineno}

\def\noi{\noindent}

\newcommand{\pd}{\gamma_P}
\newcommand{\ppt}{\operatorname{ppt}}

\newcommand{\diam}{\operatorname{diam}}

\newtheorem{thm}{Theorem}[section]
\newtheorem{cor}[thm]{Corollary}
\newtheorem{lem}[thm]{Lemma}

\theoremstyle{definition}

\theoremstyle{definition}

\theoremstyle{definition}
\newtheorem{ex}[thm]{Example}

\newcommand{\bit}{\begin{itemize}}
\newcommand{\eit}{\end{itemize}}
\newcommand{\ben}{\begin{enumerate}}
\newcommand{\een}{\end{enumerate}}
\newcommand{\beq}{\begin{equation}}
\newcommand{\eeq}{\end{equation}}
\newcommand{\bea}{\begin{eqnarray*}}
\newcommand{\eea}{\end{eqnarray*}}
\newcommand{\bpf}{\begin{proof}}
\newcommand{\epf}{\end{proof}\ms}
\newcommand{\ms}{\medskip}

\newcommand{\lc}{\left\lceil}
\newcommand{\rc}{\right\rceil}

\title{Note on power propagation time and lower bounds for the power domination number }

\author{ Daniela Ferrero\thanks{Department of Mathematics, Texas State University, San Marcos, TX, USA (dferrero@txstate.edu)}
\and Leslie Hogben\thanks{Department of Mathematics, Iowa State University, Ames, IA 50011, USA (hogben@iastate.edu) and American Institute of Mathematics, 600 E. Brokaw Road, San Jose, CA 95112, USA (hogben@aimath.org)}\and Franklin H.J. Kenter\thanks{Department of Computational and Applied Mathematics, Rice University, Houston, TX 77005, USA (franklin.kenter@gmail.com)}\and Michael Young \thanks{Department of Mathematics, Iowa State University, Ames, IA 50011, USA (myoung@iastate.edu)}}

\begin{document}


\maketitle

\begin{abstract}   We present a counterexample to a lower bound for the power domination  number given in Liao, Power domination with bounded time constraints,  {\em J. Comb. Optim.} 31 (2016): 725--742.  We also define the power propagation time, 
using the power domination propagation ideas in Liao and the (zero forcing) propagation time in Hogben et al, Propagation time for
zero forcing on a graph, {\em Discrete Appl. Math.}, 160 (2012): 1994--2005.

\end{abstract}

\noi {\bf Keywords} power domination, power propagation time, propagation time, time constraint.

\noi{\bf AMS subject classification} 05C69, 05C12, 05C15, 05C57,  94C15 


\section{Introduction}\label{sintro}

  To ensure reliability, electric power networks need to be monitored 
continuously. 
The most frequently used method of monitoring a network is to place Phase Measurement Units (PMUs) at selected electrical 
 nodes  where transmission lines, loads, and generators are connected.  A PMU placed at a node measures the voltage at the node and all current phasors at the node \cite{BMBA93}.   Because of the cost of a PMU, the trivial solution of placing a PMU at every  node is not practical, and it is
important to minimize the number of PMUs used while maintaining the
ability to observe the entire system.

This problem was first studied in terms of graphs by Haynes et al.~in
\cite{HHHH02}, where an electric power network is modeled by a
graph with the vertices representing the electric nodes and the edges
 associated with the transmission lines joining two electrical
nodes. Solving the power domination problem for a graph
consists of finding a minimum set of vertices from which the entire
graph can be observed according to rules that model PMU measurement capability; such a minimum set of  vertices will provide the locations in the
physical network where
the PMUs should be placed in order to monitor the entire network at
 minimum cost.

 Computation of the power domination number is a challenging problem. 
 Given a graph $G$ and a positive integer $k$, the problem of determining whether $G$ admits a power dominating set of cardinality at most $k$ has been proven to be NP-complete, even when restricted to bipartite graphs \cite{HHHH02}, chordal graphs \cite{HHHH02}, and planar graphs \cite{GNR08}. On the other hand, Haynes et al. gave a linear time algorithm to solve the power domination problem in trees \cite{HHHH02} and Guo et al. presented a linear time algorithm for graphs of bounded tree-width \cite{GNR08}.  

Here we reproduce the definition of the power domination process and number from \cite{L14} using our notation.  These simplified  propagation rules were shown in  \cite{BH} to be equivalent to the original propagation rules in \cite{HHHH02}. The set of neighbors of a vertex $v$ is denoted by $N(v)$, and $N(S)=\bigcup_{v\in S} N(v)$.
For a  set $S$ of vertices in a graph $G$, define  the following sets:
\ben
\item $S^{[1]}=N[S]= S\cup N(S)$.
\item For $i\ge 1$, $S^{[i+1]}=S^{[i]}\, \cup \,\{w:\exists v \in  S^{[i]},  N(v)\cap (V(G)\setminus S^{[i]})=\{w\}\}.$
\een
A set $S\subseteq V(G)$ is a {\em power dominating set} of  a graph $G$ if there is a positive integer  $\ell$ such that 
$S^{[\ell]}=V(G)$. A {\em minimum power dominating set} is a power dominating set of minimum cardinality, and the {\em power domination number}   of  $G$, denoted by $\pd (G)$, is the cardinality of a minimum power dominating set.  

For a power dominating set $S$, the {\em power propagation  time of $S$} in $G$, denoted here by $\ppt(G,S)$ and referred to as $\ell$ in \cite{AA08, L14},  is the least integer $\ell$ such that $S^{[\ell]}=V(G)$, or equivalently, the least $\ell$ such that $S^{[\ell+1]}=S^{[\ell]}$.  The {\em power propagation  time} of $G$ is 
\[\ppt(G)=\min\{\ppt(G,S)\,|\,S\mbox{ is a minimum power dominating set of } G\}. \] 
A related concept is the {\em $\ell$-round power domination number}, defined to be the minimum number of vertices needed to power-dominate $G$ in power propagation time at most $\ell$. The $\ell$-round power domination number was introduced by Aazami in \cite{Aaz, AA08}, where it was shown that its computation is an NP-hard problem even on planar graphs.

A 
lower bound on the power domination number using the propagation time of a particular power dominating set is presented in  \cite[Theorem 3]{L14}; using our notation this bound is
\beq\label{LiaoSbd}  |S|\ge \frac {|V|} { \ppt(G,S)\cdot \Delta(G)+1}
\eeq
for any power dominating set $S$ with power propagation time $\ppt(G,S)$. To make this expression useful as a lower bound for $\pd(G)$, an upper bound for $\ppt(G,S)$ must be found. It is incorrectly claimed in the proof of  \cite[Theorem 3]{L14} that $\ell=\ppt(G,S)\le\diam(G)$ for every power dominating set $S$, where  $\diam(G)$ is the diameter of $G$ (the maximum distance between any two vertices). 
Such a relationship between $\ell=\ppt(G,S)$ and $\diam(G)$ would
yield the following   incorrect lower bound for  the power domination number
\beq\label{Liaobd}  
\frac {|V|} { \diam(G)\cdot \Delta(G)+1}.\eeq

In Section \ref{serror}, we construct an infinite family of graphs having power domination number equal to 2 and having the value in \eqref{Liaobd} arbitrarily large (see Example \ref{LiaoCounterex}), thereby showing that \eqref{Liaobd} is not a lower bound for power domination number. 
 We also establish a lower bound for the power domination number of a tree that is slightly better than the one in \eqref{Liaobd}. 

We remark that there are mathematical connections between the power  domination number and the zero forcing number defined in \cite{AIM, BG} (called graph infection in the latter), and between the power propagation time discussed here and the propagation time defined in \cite{proptime}. We refer the reader to   \cite{REUF2015} for a discussion of the relationship between the power  domination number and the zero forcing number.


\section{Lower bounds for the power domination number and properties of power propagation time}\label{serror}

 
Unfortunately there is an error  in Theorem 3 in  \cite{L14}, which states:

\noi{\em  Given a connected graph $G = (V, E)$ and a positive integer time constraint $\ell$, we can derive the lower bound of the minimum cardinality of a PDS
in $G$.
$ \pd(G)\ge \frac {|V|} { \ell\cdot \Delta+1}\ge  \frac {|V|} { \diam_G\cdot \Delta+1},$
where $ \Delta$ is the maximum degree of $G$ and $\diam_G$ is the diameter of $G$.}

The first of these two inequalities is correct  if $\ell$ is the power propagation time of a \underline {minimum} power dominating set. 
Using our power propagation time terminology,
  the proof of \cite[Theorem 3]{L14} shows that
\[  |S|\ge \frac {|V(G)|} { \ppt(G,S)\cdot \Delta(G)+1}.\]
By choosing a minimum power dominating set having minimum propagating time, we  optimize this inequality (again restated in our notation):
\begin{thm}\label{LiaoOK} {\rm \cite[Theorem 3]{L14}} For a connected graph $G = (V, E)$,
\[ \pd(G)\ge \frac {|V|} { \ppt(G)\cdot \Delta(G)+1}.\]
\end{thm}

The error in the proof occurs in trying to compare the non-comparable parameters $\ppt(G)$ (or $\ppt(G,S)$) and $\diam(G)$.
The last sentence of the proof (on p. 15) reads:
 {\it Because $\ell$ is not larger than the diameter of the given graph $G$, ... }.  
Although true for trees as shown in Theorem \ref{treediam} below, this statement is not correct for an arbitrary graph, and its use yields the incorrect lower bound $ \pd(G)\ge \frac {|V(G)|} { \diam(G)\cdot \Delta+1},$
given in  \cite[Theorem 3]{L14}, to which we present a counterexample.

\begin{ex}\label{LiaoCounterex}
Given an integer $\Delta\ge 3$, construct the graph $H_\Delta$ with three levels of vertices as follows:
\ben
\item The first vertex is numbered 1.  This one vertex is on  level 1.
\item Vertex 1 has $\Delta$ neighbors, numbered $2,\dots,\Delta+1$. These vertices are on level 2.
\item Each level 2 vertex has $\Delta-1$ neighbors on level 3.  This adds $\Delta(\Delta-1)=\Delta^2-\Delta$ vertices, numbered $\Delta+2,\dots,\Delta^2+1$.
\item Add edges to make a path along level 3, so vertex $i$ is adjacent to vertex $i+1$ for $\Delta+2\le i \le \Delta^2$.
\een
Figure \ref{f:LiaoCounterex} shows the graph $H_\Delta$ for $\Delta=9$ (the three levels are the three rows of vertices).

\begin{figure}[!ht]
\begin{center}
\scalebox{.32}{\includegraphics{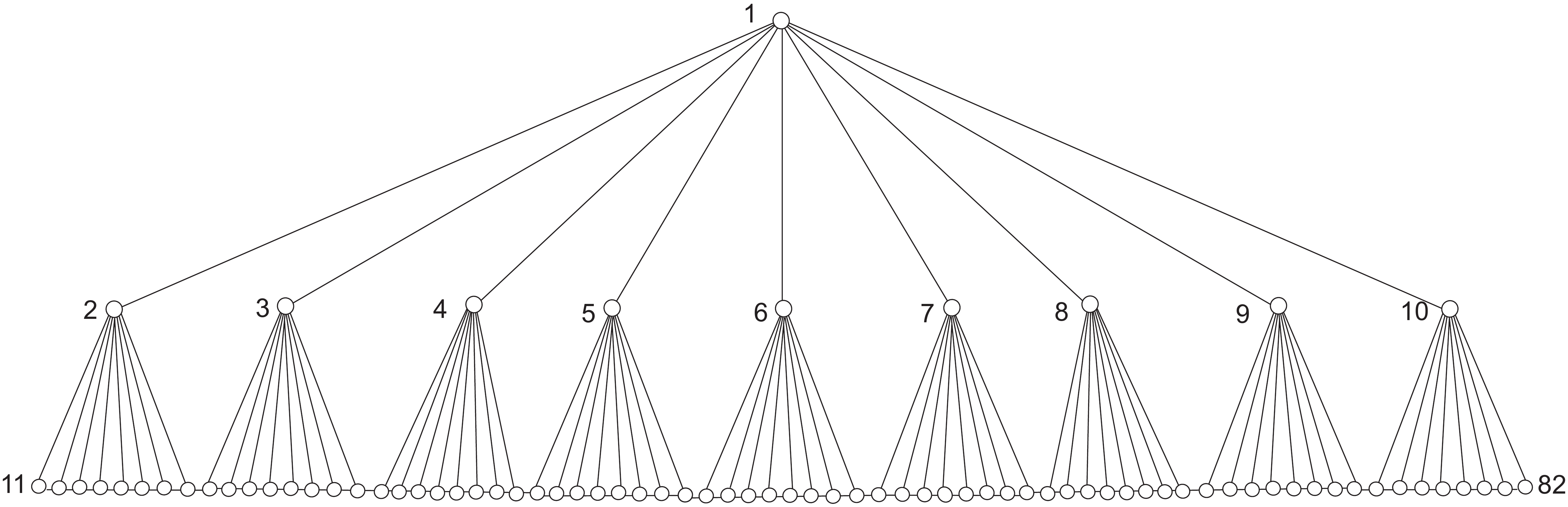}} \vspace{-8pt}
\caption{The graph $H_\Delta$ for $\Delta=9$} \label{f:LiaoCounterex}\vspace{-10pt}
 \end{center}
 \end{figure}
\enlargethispage{10 pt}

Then $|V(H_\Delta)|=1+\Delta+\Delta(\Delta-1)=1+\Delta^2,$ $\diam({H_\Delta})=4,$ and the maximum degree is $ \Delta(H_\Delta)=\Delta$.
The power domination number is $\pd(H_\Delta)=2$, because no one vertex power dominates $H_\Delta$, and vertices 1 and $\Delta+2$ are a power dominating set.
But 
\[  \frac {|V(H_\Delta)|} {\Delta(H_\Delta) \diam(H_\Delta)+1}=\frac{\Delta^2+1}{4\Delta+1}\approx \frac \Delta 4.\] By choosing $\Delta$ sufficiently large, $\frac {|V(H_\Delta)|} {\Delta(H_\Delta) \diam(H_\Delta)+1}$ can be made arbitrarily large.  In particular, for $\Delta=9$,
\[\frac {|V(H_\Delta)|} {\Delta(H_\Delta) \diam(H_\Delta)+1}=\frac {82}{37}
> 2=\pd(H_\Delta).\]
\end{ex}


While the relationship in Theorem \ref{LiaoOK} is correct, $\frac {|V(G)|} { \ppt(G)\cdot \Delta(G)+1}$ is not useful as a lower bound for $\pd(G)$, because one must know $\pd(G)$ in order to compute  $\ppt(G)$.   One can, however, use the relationship in Theorem \ref{LiaoOK}  to obtain a lower bound on power propagation time, assuming one knows $\pd(G)$:

\begin{cor}\label{pptbd}  For a graph $G$,\vspace{3pt}
\[ \ppt(G)\ge\lc \frac {|V(G)|-\pd(G)} {\pd(G) \cdot \Delta(G)}\rc.\]
\end{cor}

  Although in general it is not true that $\ppt(G)\le \diam(G)$,  we show below that  $\ppt(T)\le \diam(T)$ when $T$ is a tree, and in fact this can be strengthened slightly when $T$ has at least three vertices 
  (see Theorem \ref{treediam}). 
First we establish a relationship between the power propagation time and the maximum  length of a trail.  
 A {\em walk} $v_0v_1\cdots v_p$ in $G$ is a subgraph with vertex set $\{v_0,v_1,\dots, v_p\}$  and edge set $\{v_0v_1, v_1v_2, \dots, v_{p-1}v_p\}$ (vertices and/or edges may be repeated in  a walk but duplicates are removed from sets). A  {\em trail} is a walk with no repeated edges (vertices may be repeated).   
Of course, a {path} is a trail with no repeated vertices. The {\em length} of a trail $P=v_0v_1\cdots v_p$ is $p$, i.e., the number of edges in $P$.  

Given a power dominating set $S$, the power propagation process naturally labels a vertex by the time at which it is observed, that is, $t(v)=i$ for $v\in S^{[i]}\setminus S^{[i-1]}$ with $i\ge 1$ and $t(v)=0$ for $i\in S$.  
We can extend this labeling to edges by  $t(uv)=\max\{ t(u),t(v)\}$.  A trail $P=v_0v_1\cdots v_p$ is called {\em monotone}  if $p\ge 1$, $t(v_{i-1}v_{i})\le t(v_{i}v_{i+1})\le t(v_{i-1}v_{i})+1$ for $i=1,\dots,p-1$ and $t(v_{p-1}v_p)=t(v_p)$.  Note that the term `monotone' comes from the monotonicity of the edge labels.

\begin{lem}\label{ztraillem}  Let $G$ be a graph    
and let $S$ be a power dominating set of $G$ such that $S$ has no vertices of degree $1$. Then for every vertex $v \in V(G)\setminus S$, there is a monotone trail of length   at least  $t(v) +1$ in which $v$ is the last  vertex of the trail. 
\end{lem}
\bpf  
Suppose $t(v) = 1$.  Then there is a vertex $u\in S\cap N(v)$, and since $S$ has no degree 1 vertices, there is another neighbor $w$ of $u$.  Since $t(u)=0$ and $t(w) = 0$ or 1, $P=wuv$ is the required trail.  This proves the base case. 

Assume the statement is true when $1\le t(v) < i$, and suppose  that $t(v) = i\ge 2$. 
 Let $w \in V(G)$ be the vertex that forces $v$ at time $i$. Assume first that $t(w)=i-1$.  Then  there is a monotone trail $P$ of length at least $i$ that has $w$ as its last vertex. Since $t(wv)=i$,  we can adjoin the edge $wv$ and vertex $v$ to $P$ to obtain  a monotone trail of length at least $i+1$ with $v$ as last vertex. 
 Now assume that $1\le t(w)<i-1$.  Then there exists a vertex $w' \in N(w)$ such that $t(w') = i-1$. Therefore, there is a monotone trail $P$ of length at least $i$ that has $w'$ as its last vertex. Since $t(w'w)=i-1$ and $t(wv)=i$,  we can adjoin edge $w'w$, vertex $w$, edge $wv$, and vertex $v$ to $P$ to obtain   a monotone trail $P'$ of length at least $i+2$ that ends at $v$. 
  \epf


\begin{thm}\label{treediam}  For every tree $T$ with at least $3$ vertices, $T$ has a minimum power dominating set $S$ with $\ppt(T,S) = \ppt(T)$ and $\deg_T u\ge 2$ for all $u\in S$.  Furthermore, $\ppt(T)\le \diam(T)-1$.
\end{thm}

\bpf
Let $S$ be a power dominating set of $T$ such that $\ppt(T,S) = \ppt(T)$. Assume $S$ contains a vertex $v$ with degree 1. Let $u$ be the vertex adjacent to $v$ in $T$. The degree of $u$ is at least 2, because $T$ has at least 3 vertices. Since $S$ is minimum, $u \notin S$. Therefore $S' = S \setminus \{v\} \cup \{u\}$ is also a minimum power dominating set of $T$, and $\ppt(T,S') \le \ppt(T,S)$. So $T$ has a minimum power dominating set containing no vertices of degree 1 that has power propagation time of $\ppt(T)$. By Lemma \ref{ztraillem}, $T$ contains a monotone trail of length at least $\ppt(T) + 1$, and  every trail in  a tree is a path. Therefore, $T$ has a path of length at least $\ppt(T) + 1$. So $\ppt(T) + 1 \le \diam(T)$.
\epf

\begin{cor}\label{treediamdpd}  For every tree $T$ with at least $3$ vertices, 
\[ \pd(T)\ge \lc\frac {|V(G)|} { (\diam(G)-1)\cdot \Delta(G)+1}\rc.\]
\end{cor}




\end{document}